\documentclass[11pt]{article}
\usepackage{amsmath}
\usepackage{amstext,amsbsy}
\usepackage{epsfig,multicol}
\usepackage{amsthm}
\usepackage{amssymb,latexsym}
\usepackage{color}
\usepackage{tikz}
\usepackage{subfigure}
\usepackage{url}
\usepackage{enumerate}
\usetikzlibrary{calc}
\usetikzlibrary{arrows}

\newtheorem{thm}{Theorem}%[section]
\newtheorem{defi}{Definition}
\newtheorem{exa}[thm]{Example}
\newtheorem{lemma}[thm]{Lemma}
\newtheorem{cor}[thm]{Corollary}
\newtheorem{prop}[thm]{Proposition}

\newcommand{\match}{\mathcal{M}}

\newcommand{\lmark}{\mathcal{L}}
\newcommand{\bmark}{\mathcal{B}}
\newcommand{\lmm}{\ell}
\newcommand{\bmm}{b}

\title{Counting with Borel's Triangle}
\author{Yue Cai\thanks{ycai@math.tamu.edu} \ and \   Catherine Yan\thanks{cyan@math.tamu.edu}\\
Department of Mathematics, Texas A\&M University, College Station, TX 77843}
\date{}

\begin{document}

\maketitle 

\begin{abstract} 
Borel's triangle is an array of integers closely related to the classical Catalan numbers. 
In this paper we study combinatorial statistics counted by  Borel's triangle.  
We present various combinatorial interpretations of Borel's triangle in terms of lattice paths, binary trees, 
and pattern avoiding permutations and matchings,  and derive a functional equation that 
is useful in analyzing the involved structures.  
\end{abstract}



\section{Introduction}
%[Introducing the background of this two sequences.  History, our motivation and outline of results. 
%Organization of the paper ]

In combinatorics, Catalan's triangle is a triangular array whose right boundary, as well as the row-sum, 
are the classical Catalan sequence. 
%It can also be viewed as a refinement of the Catalan numbers since the row-sum of this array is also  the Catalan number.  
Entries in  Catalan's  triangle appear in many combinatorial structures, notably in lattice paths, 
plane trees and binary trees, triangulations, and parking functions. 
See the references at the On-line Encyclopedia of Integer Sequences (OEIS) \cite{OEIS} for the 
sequences A009766 and A033184, which are   Catalan's  triangle and its transpose, respectively. 

Recently another triangular array that is closely related to Catalan's triangle 
has  appeared in various studies in commutative algebra, combinatorics,   and discrete 
geometry.  It is the sequence A234950 in OEIS  and is called \emph{Borel's triangle}, which is 
related to pseudo-triangulations of point sets \cite{AOSS08} and the Betti numbers of certain principal Borel ideals \cite{FMS11},  and appears in Cambrian Hopf algebras \cite{Chatel-Pilaud},  quantum physics \cite{LPZ14}, 
and  permutation patterns \cite{Remmel14}.  
  %The name ``Borel's triangle"  was first appeared in \cite{FMS15} by Francisco, Martin and Schweig since this 
 %array represents the Betti numbers of the smallest Borel ideals of $K[x_1, x_2, \dots, x_n]$ containing 
% the monomial  $x_1x_2x_3\cdots x_n$.  
In the second author's work of  parking functions and parking distributions on trees, Borel's triangle gives   
the coefficients of certain generating functions on the nondecreasing parking functions \cite[Section 3]{BGY17},
which inspires the project on finding classes of objects  that are counted by  Borel's triangle and characterizing their combinatorial structures.  

In this paper we  study Borel's triangle from combinatorial and enumerative contexts. 
We start by  recalling the  preliminary results of Catalan's and Borel's triangles in Section 2. 
Then we describe combinatorial interpretations of Borel's triangle in terms of marked Catalan structures, 
vertex-marked binary trees, and combinatorial statistics of  matchings avoiding certain pairs of patterns. 
This is the content of Section 3. 
In Section 4 we present a bijection between restricted Dumont permutations and marked Dyck paths; the latter is a typical structure enumerated by Borel's triangle. This bijection leads to a functional equation for Borel's triangle, which is used in  the last section to analyze other structures counted by Borel's triangle.

% Preliminary results about Catalan's and Borel's triangles. 
\section{Preliminary Results on Catalan's and Borel's Triangles} 
\subsection{Catalan's Triangle} 
\emph{Catalan's triangle} $\{C_{n,k}:\ 0 \leq k \leq n\}$  
is the array defined 
by the recurrence $C_{n,k}= C_{n-1,k}+ C_{n,k-1}$ for $0 < k < n$ and the boundary conditions 
$C_{n,0}=1$ and $C_{n,n}=C_{n,n-1}$.  The first seven rows of Catalan's triangle are given below. 

\begin{center} 
\begin{tabular}{l|ccccccc   } 
 $n\setminus k $ & 0 & 1 & 2 & 3 & 4 & 5  & 6 \\  \hline 
 0               & 1 &   &   &   &   &    &   \\ 
 1               & 1 & 1 &   &   &   &    &   \\
 2               & 1 & 2 & 2 &   &   &    &   \\
 3               & 1 & 3 & 5 & 5 &   &    &    \\
 4               & 1 & 4 & 9 & 14& 14&    &    \\
 5               & 1 & 5 & 14& 28& 42& 42 &    \\
 6               & 1 & 6 & 20 & 48 & 90 & 132 & 132 
 \end{tabular}
\end{center}

The entries in Catalan's  triangle are often called ballot numbers, since $C_{n,k}$ counts the number 
of lattice paths in the coordinate plane from $(0,0)$ to $(n,k)$ that do not go above the line $y=x$. Explicitly, 
\[
 C_{n,k} = \frac{n-k+1}{n+1} \binom{n+k}{n}. 
\]
In particular,  $C_{n,n}=\frac{1}{n+1} \binom{2n}{n}$ is the $n$-th Catalan number.  
The sum of entries in the $n$-th row is the $(n+1)$-st Catalan number. 
The bivariate generating function $\mathcal{C}(t,x) = \sum_{n,k} C_{n,k} t^k x^n$ 
can be expressed as 
\begin{equation} \label{eq1} 
\mathcal{C}(t,x)= \frac{C(tx)}{1-xC(tx)}, 
\end{equation}
where $C(x)$ is the generating function for Catalan numbers, i.e., 
\begin{equation} \label{Catalan-gf} 
 C(x) = \sum_{n=0}^\infty \frac{1}{n+1} \binom{2n}{n} x^n = \frac{1-\sqrt{1-4x}}{2x}. 
\end{equation}

Catalan's  triangle appears in countless places throughout enumerative combinatorics. 
A few examples are given in Theorem \ref{Catalan-array}, which can be found directly or derived 
from the examples in OEIS A009766. 
As usual a \emph{Dyck path of semi-length} 
$n$ is a  lattice path in the coordinate plane from  $(0, 0)$ to $(2n, 0)$ consisting of $n$ up-steps (along the
vector $(1, 1)$)  and $n$ down-steps (along the vector $ (1, -1)$)  such that the path never goes 
below the $x$-axis.  A \emph{parking function of length $n$}  is a sequence of positive integers $(a_1, \dots, a_n)$ such that $a_i \leq \sigma_i$ for some permutation $\sigma \in \mathfrak{S}_n$. 

\begin{thm} \label{Catalan-array} 
The entry $C_{n,k}$ of  Catalan's  triangle counts  the following sets. The configurations with $n=3$ and $k=2$ are listed. 
\begin{enumerate}[(1)] 
\item Lattice paths in the coordinate plane from $(0,0)$ to $(n,k)$ that never go above the line $y=x$.
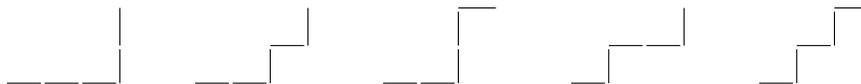
\begin{center} 
\begin{tikzpicture}[scale=0.5]
\draw (0,0)--(0.9,0) (1,0)--(1.9,0) (2,0)--(2.9,0) (3,0)--(3,0.9) (3,1)--(3,2) 
(5,0)--(5.9,0) (6,0)--(6.9,0) (7,0)--(7,0.9) (7,1)--(7.9,1) (8,1)--(8,2) 
(10,0)--(10.9,0) (11,0)--(11.9,0) (12,0)--(12,0.9) (12,1)--(12,1.9) (12,2)--(13,2)
 (15,0)--(15.9,0) (16,0)--(16,0.9) (16,1)--(16.9,1) (17,1)--(17.9,1) (18,1)--(18,2) 
(20,0)--(20.9,0) (21,0)--(21,0.9) (21,1)--(21.9,1) (22,1)--(22,1.9) (22,2)--(23,2); 
\end{tikzpicture}
\end{center} 

\item  Dyck paths of semi-length $n+1$ that have $k$ up-steps (or down-steps)  not at ground level.  Equivalently,
it is the set of Dyck paths of semi-length $n+1$ with $n+1-k$ returns to the $x$-axis, $($not counting the 
starting point $(0,0)$.$) $
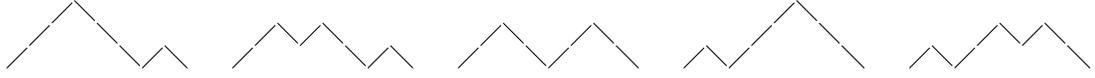
\begin{center} 
\begin{tikzpicture}[scale=0.3] 
\draw (0,0)--(0.9,0.9) (1,1)--(1.9,1.9) (2,2)--(2.9,2.9) (3,3)--(3.9,2.1) (4,2)--(4.9, 1.1) (5,1)--(5.9,0.1)(6,0)--(6.9,0.9) (7,1)--(8,0) 
(10,0)--(10.9,0.9) (11,1)--(11.9,1.9) (12,2)--(12.9, 1.1) (13,1)--(13.9, 1.9) (14,2)--(14.9,1.1) (15,1)--(15.9,0.1) (16,0)--(16.9,0.9) (17,1)--(18,0) 
(20,0)--(20.9,0.9) (21,1)--(21.9, 1.9) (22,2)--(22.9,1.1) (23,1)--(23.9,0.1) (24,0)--(24.9,0.9) (25,1)--(25.9,1.9) (26,2)--(26.9,1.1) (27,1)--(28,0) 
(30,0)--(30.9,0.9) (31,1)--(31.9,0.1) (32,0)--(32.9,0.9) (33,1)--(33.9,1.9) (34,2)--(34.9, 2.9) (35,3)--(35.9, 2.1) (36,2)--(36.9, 1.1) (37,1)--(38,0) 
(40,0)--(40.9,0.9) (41,1)--(41.9,0.1) (42,0)--(42.9,0.9) (43,1)--(43.9, 1.9) (44,2)--(44.9, 1.1) (45,1)--(45.9,1.9) (46,2)--(46.9,1.1) (47,1)--(48,0); 
\end{tikzpicture}
\end{center} 

\item Dyck paths of semi-length $n+1$ and having the first (or the last) peak at height $n-k+1$. 
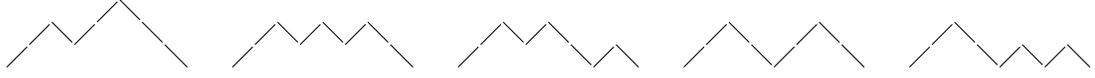
\begin{center} 
\begin{tikzpicture}[scale=0.3] 
\draw  (0,0)--(0.9,0.9) (1,1)--(1.9,1.9) (2,2)--(2.9,1.1) (3,1)--(3.9,1.9) (4,2)--(4.9,2.9) (5,3)--(5.9,2.1) (6,2)--(6.9, 1.1) (7,1)--(8,0)
 (10,0)--(10.9,0.9) (11,1)--(11.9,1.9) (12,2)--(12.9,1.1) (13,1)--(13.9,1.9) (14,2)--(14.9,1.1) (15,1)--(15.9,1.9) (16,2)--(16.9,1.1) (17,1)--(18,0) 
(20,0)--(20.9,0.9) (21,1)--(21.9,1.9) (22,2)--(22.9,1.1) (23,1)--(23.9, 1.9) (24,2)--(24.9,1.1) (25,1)--(25.9,0.1) (26,0)--(26.9,0.9) (27,1)--(28,0) 
(30,0)--(30.9,0.9) (31,1)--(31.9,1.9) (32,2)--(32.9,1.1) (33,1)--(33.9, 0.1) (34,0)--(34.9,0.9) (35,1)--(35.9,1.9) (36,2)--(36.9, 1.1) (37,1)--(38,0)
(40,0)--(40.9,0.9) (41,1)--(41.9, 1.9) (42,2)--(42.9,1.1) (43,1)--(43.9, 0.1) (44,0)--(44.9, 0.9) (45,1)--(45.9, 0.1) (46,0)--(46.9, 0.9) (47,1)--(48,0); 
\end{tikzpicture}
\end{center} 

\item Unlabeled plane  forests on $n+1$ vertices such that there are  $n+1-k$ planted plane trees. 
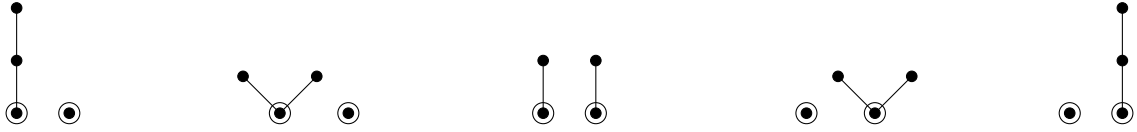
\begin{center} 
\begin{tikzpicture}[scale=0.7] 
\draw [fill] (0,0) circle [radius=0.1]; \draw (0,0) circle [radius=0.2]; 
\draw [fill] (1,0) circle [radius=0.1]; \draw (1,0) circle [radius=0.2]; 
\draw [fill] (0,1) circle [radius=0.1];
\draw [fill] (0,2) circle [radius=0.1];
\draw (0,0)--(0,2); 

\draw [fill] (5,0) circle [radius=0.1]; \draw (5,0) circle [radius=0.2]; 
\draw [fill] (4.3,0.7) circle [radius=0.1];
\draw [fill] (5.7,0.7) circle [radius=0.1];
\draw [fill] (6.3,0) circle [radius=0.1];  \draw (6.3,0) circle [radius=0.2]; 
\draw (4.3,0.7)--(5,0)--(5.7,0.7); 

\draw [fill] (10,0) circle [radius=0.1]; \draw (10,0) circle [radius=0.2]; 
\draw [fill] (11,0) circle [radius=0.1]; \draw (11,0) circle [radius=0.2]; 
\draw [fill] (10,1) circle [radius=0.1];
\draw [fill] (11,1) circle [radius=0.1];
\draw (10,0)--(10,1) (11,0)--(11,1);

\draw [fill] (15,0) circle [radius=0.1]; \draw (15,0) circle [radius=0.2]; 
\draw [fill] (16.3,0) circle [radius=0.1]; \draw (16.3,0) circle [radius=0.2]; 
\draw [fill] (15.6,0.7) circle [radius=0.1];
\draw [fill] (17,0.7) circle [radius=0.1];
\draw (15.6,0.7)--(16.3,0)--(17, 0.7); 

\draw [fill] (20,0) circle [radius=0.1]; \draw (20,0) circle [radius=0.2]; 
\draw [fill] (21,0) circle [radius=0.1]; \draw (21,0) circle [radius=0.2]; 
\draw [fill] (21,1) circle [radius=0.1];
\draw [fill] (21,2) circle [radius=0.1];
\draw (21,0)--(21,2); 
\end{tikzpicture}
\end{center} 

\item Nondecreasing parking functions of length $n+1$ with maximal element $k+1$. 
\[ 1113 \qquad 1123 \qquad 1133 \qquad 1223 \qquad 1233 \]

\item Nondecreasing parking functions of length $n+1$ with $k$ unlucky ``cars", 
 that is, the entries $a_i$  such that  $a_i \neq i$. \footnote{In an ordinary parking 
 function $(a_1\dots, a_n)$ defined via a parking process \cite{Yan14},   a car $C_i$ is \emph{lucky}  if  $C_i$ is parked at its preferred space. For nondecreasing parking functions, the car $C_i$ is lucky if and only if $a_i=i$.  } 
\[ 1114 \qquad 1124 \qquad 1133 \qquad 1222 \qquad 1223 \]

\item Nondecreasing parking functions of length $n+1$ containing $n-k+1$ ones. 
\[ 1122 \qquad 1123 \qquad 1124 \qquad 1133 \qquad 1134\]
\end{enumerate} 
\end{thm}

\subsection{Borel's Triangle} 
%$B_{n,k}$, relation with $C_{n,k}$, G.F., and several known meanings. 

{Borel's  triangle} $\{B_{n,k}:\  0 \leq k \leq n\}$ 
is an array of numbers obtained from an invertible  transformation to Catalan's triangle 
by the equation 
\begin{equation} \label{eq-Borel} 
B_{n,k}= \sum_{s=k}^n \binom{s}{k} C_{n,s}. 
\end{equation} 
 %where $C_{n,s}$ is entries in Catalan's  triangle. 
Equivalently, 
\[
\sum_{k=0}^n B_{n,k} t^k = \sum_{k=0}^n C_{n,k} (1+t)^k, 
\]
and the bivariate generating function 
$\mathcal{B}(t,x):= \sum_{n, k \geq 0} B_{n,k} t^k x^n $ satisfies 
\begin{equation} \label{eq-Borel2} 
\mathcal{B}(t,x) = \mathcal{C}(1+t, x)
 = \frac{C((1+t)x)}{1-xC((1+t)x)},
\end{equation} 
where $C(x)$ is the Catalan generating function given in \eqref{Catalan-gf}. 
A variation of Borel's triangle,   the OEIS sequence A094385,  appeared in Barry's study of generalized Pascal
 matrices defined by Riordan's arrays \cite{Barry13}.  Barry gave 
 an explicit expression of $B_{n,k}$ as
 \begin{equation} \label{eq2.2-1} 
 B_{n,k} = \frac{1}{n+1} \binom{2n+2}{n-k} \binom{n+k}{n}. 
 \end{equation} 
The entries of $B_{n,k}$ for small values of $n$ and $k$ are listed below. 

\begin{center} 
\begin{tabular}{l|ccccccc   } 
 $n\setminus k $ & 0 & 1 & 2 & 3 & 4 & 5  & 6 \\  \hline 
 0               & 1 &   &   &   &   &    &   \\ 
 1               & 2 & 1 &   &   &   &    &   \\
 2               & 5 & 6 & 2 &   &   &    &   \\
 3               & 14 & 28 & 20 & 5 &   &    &    \\
 4               & 42 & 120 & 135 & 70& 14&    &    \\
 5               & 132 & 495 & 770 & 616 & 252 & 42 &    \\
 6               & 429 & 2002  & 4004  & 4368  & 2730  & 924  & 132 
 \end{tabular}
 \end{center}

Computing the   sum of the entries in each row of  Borel's triangle  and adding an extra $1$ at the beginning, 
we obtain the sequence $1, 1, 3, 13, 67, 381, 2307, 14589, \dots$. This sequence is 
called the \emph{generalized Catalan Number} and is denoted by  $C(2,n)$ in OEIS (A064062).  
Explicitly,  $C(2,n)$ is the sum of the entries in the $(n-1)$-st row of Borel's triangle, i.e., 
\begin{equation} \label{eq2.2-2} 
C(2,i)=\sum_{k} B_{i-1, k} = \sum_{k} C_{i-1, k}  2^k, 
\end{equation} 
and 
\begin{equation}  \label{eq2.2-3} 
\sum_{i\geq 0} C(2,i)  x^i  = \frac{1+2x C(2x)}{1+x} = \frac{1}{1-xC(2x)} 
=\frac{4}{3+ \sqrt{1-8x}}. 
\end{equation}

\section{Combinatorial Interpretations of Borel's Triangle} 
\label{sec: Interpretation}

In this section we present various interpretations of Borel's triangle and the generalized 
Catalan numbers in basic combinatorial structures,  in particular, in Catalan structures, 
binary trees, and permutations. 

\subsection{Marked Catalan Structures}   

%\input{1marked-Catalan.tex}

%input file for subsection "Marked Catalan Structure"

First we note that  using the defining relation \eqref{eq-Borel} , we can easily get a family of  
combinatorial interpretations of $B_{n,k}$ from those of $C_{n,k}$. 
Assume that we have a discrete structure $P$ with a combinatorial statistic $\alpha: P \rightarrow \mathbb{N}$ 
counting some special features of $P$. 
If $C_{n,k}$ counts the set of such structures with  parameter $n$ (usually indicates the size)  for which 
$\alpha(P)=k$, then the generalized Catalan number $C(2,n)$  counts the set of pairs $\{(P, S)\}$ in which 
 $P$ is a structures with parameter $n$ and $S$ is a subset of the special features of $P$, and 
 $B_{n,k}$ counts those pairs  $\{(P, S)\}$ with $|S|=k$. 
We say that the special features in $S$ are \emph{marked}, and the set of pairs $(P, S)$ 
a \emph{marked Catalan structure},
which are the most basic structures counted by $B_{n,k}$.   
We will establish other combinatorial interpretations of $B_{n,k}$ by constructing bijections to a 
marked Catalan structure. 

Each item in Theorem \ref{Catalan-array} gives a marked Catalan structure. 
In particular,  in this paper  we would frequently use the following one, which comes from item (2) 
of Theorem \ref{Catalan-array}.

%[to be checked . which is used ? ] 
\begin{thm} \label{Theorem:Borel-marks}
The entry $B_{n,k}$ of Borel's triangle counts the set of  $(D, S)$ where $D$ is a 
Dyck path of semi-length $n+1$ and $S$ consists of  $k$  up-steps  of $D$,  none of which is  at ground level. 
 Consequently, the generalized Catalan number $C(2,n)$ counts all the pairs $(D, S)$  where $D$ is a 
Dyck path of semi-length $n+1$ and $S$ is a set of up-steps of $D$ not at ground level. 
\end{thm} 
%We say that the up-steps in $S$ are marked.
By symmetry of Dyck paths, the ``up-steps" in the previous statement can be replaced by ``down-steps".

\subsection{Binary trees with Marked Vertices}

%\input{2binary-tree.tex} 

% binary.-tree.tex 

In \cite{FMS15} Francisco et al.\ gave four  interpretations of the entries of Borel's triangle; 
two of them are  in terms of marked binary trees, which are shown to be counted by $B_{n,k}$ via 
 bijections to two structures in commutative algebra and  discrete geometry, namely, the EK-symbols of the smallest 
 Borel ideal of $K[x_1, x_2, \dots, x_n ] $ containing the monomial $x_1x_2\cdots x_n$, and the pointed 
 pseudotriangulations, 
 The two  interpretations of marked binary trees  are different from 
the marked Catalan structures  in the preceding subsection  as  each set contains binary trees with $k$ marked 
vertices and $n+1$ unmarked vertices,  and hence the total number of vertices in the binary trees is  $n+k+1$. 
Our goal is to find more direct, combinatorial  proofs for these two sets of 
marked binary trees.

First we recall  the definition of marked binary trees from \cite{FMS15}. 
A \emph{binary tree} is a rooted plane tree in which each vertex has at most two children, 
which are referred to as the \emph{left child} and the \emph{ right child}.  A vertex of a binary 
tree is a \emph{leaf} if it is has  no children, and it is  a \emph{branching vertex} if it has two
children.  

%Following \cite{FMS15}, we define the leaf-marked and branch-marked binary trees as follows. 

\begin{defi}\emph{
Let $T$ be a binary tree. The \emph{rightmost} leaf of $T$ is the last leaf of $T$ in the preorder depth-first traversal,  i.e., the  leaf obtained by starting at the 
root and descending, taking a right edge whenever possible.  Let $X$ be a set of leaves of $T$ not containing the rightmost leaf. Then the pair $(T,X)$ is called a \emph{leaf-marked binary tree}.  We say that the 
vertices in $X$ are \emph{marked}, and the other vertices \emph{unmarked}. }

\emph{If $B$ is a set of branching vertices of $T$, we call the pair $(T, B)$ a \emph{branch-marked 
binary tree}. } 
\end{defi} 

%Fransisco et al. proved the following two combinatorial interpretations for $B_{n,k}$. 
\begin{thm}\cite{FMS15}   \label{thm:FMS} 
The entry $B_{n,k}$ of Borel's triangle counts the following sets.    
\begin{enumerate}[(a)] 
\item  Leaf-marked binary trees with $n+1$ unmarked vertices and $k$ marked leaves.
%where the last leaf in the pre-order depth-first search traversal is never marked. 
\item Branch-marked binary trees with $n+1$ unmarked vertices and $k$  marked 
branching vertices. 
\end{enumerate} 
\end{thm} 
%Note that in each case of (a) or (b), the total number of vertices in the binary trees is  $n+k+1$. 
As a corollary, the generalized Catalan number $C(2,n)$
 counts the number of leaf-marked binary trees with $n$  unmarked vertices or branch-marked binary trees with 
 $n$ unmarked vertices.  (There is no constraint on the number of marked vertices.)

In this section we  give a bijective proof for (a).   A proof of (b) is presented in Section \ref{Sec:functional} 
using a tree decomposition and a functional equation of Borel's triangle. We remark that 
an easy  bijection between the sets in  (a) and (b) is given in \cite{FMS15},
which is based on the observation that any binary tree with $k$ branching vertices must have $k+1$ leaves,
hence one can associate the branching vertices bijectively to the set of all leaves except the rightmost one. 
This bijection, although very simple,  does not give much insight on  the structure of 
branching vertices.   

Let $U$ represents an up-step and $D$ represents a down-step in a Dyck path. A Dyck path $P$ of semi-length $n$  is represented as a $\{U,D\}$-sequence of length $2n$ with 
exactly $n$  $U$'s and $n$ $D$'s. Any such a sequence must end with a sequence of down-steps, which 
we call the \emph{last run} of down steps of $P$. 
A \emph{$UDD$-pattern} is a subword of $P$ consists of three consecutive letters $U$, $D$, $D$.  Let $S$ be a set 
of $UDD$-patterns of a Dyck path $P$ in which no pattern contains any down-steps  in the last run of down-steps of
$P$. 
We say that the pair $(P, S)$ is a \emph{Dyck path with marked $UDD$-patterns}.
We shall  prove that Dyck paths with marked $UDD$-patterns are also counted by Borel's triangle, and then we 
construct a bijection between  the set of leaf-marked 
binary trees and  that of Dyck paths with marked $UDD$-patterns.

\begin{prop} \label{UDD-dyck} 
The entry $B_{n,k}$ of Borel's triangle counts the number of Dyck paths  of semi-length $n+k+1$
with $k$ marked $UDD$-patterns. 
\end{prop} 
\begin{proof} 
We show that the number of Dyck paths  of semi-length $n+k+1$ with $k$ marked $UDD$-patterns
can be computed by the formula $\sum_{s=k}^n \binom{s}{k} C_{n,s}$. Then the claim follows from 
Eq.~\eqref{eq-Borel}. 

By item (3) of Theorem 1, $C_{n,s}$ is the number of Dyck paths of semi-length $n+1$ with the last
peak at height $n+1-s$. Any such a Dyck path has $s$ 
down-steps not in the last run.  If $s \geq k$, we can pick $k$ down-steps not in the last run and 
replace each of them with a marked $UDD$. This results in a Dyck path of semi-length $n+k+1$ with 
$k$ marked $UDD$-patterns. Conversely, from any Dyck path of semi-length $n+k+1$ with $k$ marked $UDD$s 
that are disjoint from the last run of down steps, replacing each marked $UDD$ with $D$ would yield a Dyck 
path of semi-length  $n+1$
with at least $k$ down steps not in the last run.   
\end{proof} 

\begin{thm} 
There is a bijection between the set of leaf-marked binary trees with $n+1$ unmarked vertices 
and $k$ marked leaves  and the set of Dyck paths of semi-length $n+k+1$ with $k$ marked $UDD$-patterns. 
\end{thm} 
\begin{proof} 
A full binary tree is a binary tree in which every vertex has either 0 or 2 children.  
Given a binary tree $T$ with $n+k+1$ vertices, first we turn it into a full binary tree 
$\mathcal{T}$ by adding left and/or right children to any non-branching vertex so that every 
vertex of $T$ has two children.  Any mark in $T$ remains unchanged. 
In the full binary tree $\mathcal{T}$ the vertices with two children 
are called \emph{internal vertices}, and 
$\mathcal{T}$ contains $n+k+1$ internal vertices and $n+2+k$ leaves.  Clearly one can recover $T$ 
from $\mathcal{T}$ by removing all the leaves. 

Both the set of full binary trees with $m$ interval vertices and the set of Dyck paths of semi-length 
$m$ are counted by the $m$-th Catalan number, and there are well-known bijections between these two sets. 
One bijection goes as follows. Given  a full binary tree $\mathcal{T}$,  we start at the root and traverse  in the preorder 
traversal. In the process  record ``$U$" for each internal  vertex  and ``$D$" for each leaf of $\mathcal{T}$, 
except for the rightmost leaf of $\mathcal{T}$.  This gives a $\{U, D\}$-sequence that corresponds to a 
Dyck path of semi-length $n+k+1$, see \cite[Corollary 6.2.3]{EC2}.

 A leaf of $T$  corresponds to an internal vertex of $\mathcal{T}$ for which both children in 
 $\mathcal{T}$ are leaves. Under the above bijection it  corresponds to a $UDD$-pattern in the word 
 representation of the Dyck path except for the rightmost leaf of $T$ which corresponds to the last 
 peak of the Dyck path. Hence if $S$ is a set of leaves of $T$ not containing the rightmost one, then 
 under the above map it becomes a set of $UDD$-patterns not containing any down steps in the last run. 
\end{proof}

\begin{figure}[ht] 
\begin{center} 
\begin{tikzpicture} [scale=0.7] 
\draw (2.5, 4.5)--(1,3)--(0,2) (2.5,4.5)--(4,3)--(5,2)--(4,1) (1,3)--(2,2)--(1,1) (2,2)--(3,1); 
\draw [fill] (2.5,4.5) circle [radius=0.1];
\draw [fill] (1,3) circle [radius=0.1];
\draw [fill=red,red] (0,2) circle [radius=0.1]; \draw [red] (0,2) circle [radius=0.2]; 
\draw [fill] (4,3) circle [radius=0.1];
\draw [fill] (5,2) circle [radius=0.1];
\draw [fill] (4,1) circle [radius=0.1];
\draw [fill] (2,2) circle [radius=0.1];
\draw [fill] (1,1) circle [radius=0.1];
\draw [fill=red,red] (3,1) circle [radius=0.1]; \draw [red] (3,1) circle [radius=0.2]; 
\node at (2.5,-0.5) {a binary tree $T$ with marked leaves} ;

\draw (12.5, 4.5)--(11,3)--(10,2) (12.5,4.5)--(14,3)--(15,2)--(14,1) (11,3)--(12,2)--(11,1) (12,2)--(13,1); 
\draw [fill] (12.5,4.5) circle [radius=0.1];
\draw [fill] (11,3) circle [radius=0.1];
\draw [fill=red,red] (10,2) circle [radius=0.1]; \draw [red] (10,2) circle [radius=0.2]; 
\draw (10,2)--(9.3,1.3) (10,2)--(10.7, 1.3); 

\draw [fill] (14,3) circle [radius=0.1]; \draw (14,3)--(13.3, 2.3); 
\draw [fill] (15,2) circle [radius=0.1]; \draw (15,2)--(15.7, 1.3); 

\draw [fill] (14.2,1.2) circle [radius=0.1]; \draw (14.2,1.2)--(13.6,0.6) (14.2,1.2)--(14.9,0.5) ;
\draw [fill] (12,2) circle [radius=0.1];
\draw [fill] (11,1) circle [radius=0.1]; \draw (11,1)--(10.3,0.3) (11,1)--(11.7,0.3) ;
\draw [fill=red,red] (13,1) circle [radius=0.1]; \draw [red] (13,1) circle [radius=0.2]; 
\draw (13,1)--(12.3,0.3) (13,1)--(13.7,0.3) ;
\node at (13.5,-0.5) {the  full  binary tree $\mathcal{T}$,} ;

\draw [scale=0.8] 
(4,-5)--(4.9,-4.1) (5,-4)--(5.9,-3.1) (6,-3)--(6.9, -2.1) (7,-2)--(7.9, -2.9) (8,-3)--(8.9, -3.9)
(9,-4)--(9.9, -3.1)
 (10,-3)--(10.9,-2.1) (11,-2)--(11.9,-2.9) (12,-3)--(12.9,-3.9) (13,-4)--(13.9, -3.1)
 (14,-3)--(14.9, -3.9) (15,-4)--(15.9,-4.9) 
(16,-5)--(16.9, -4.1) (17,-4)--(17.9, -4.9) (18,-5)--(18.9, -4.1) 
(19,-4)--(19.9, -3.1) (20,-3)--(20.9, -3.9) (21,-4)--(22,-5); 

\draw [scale=0.8, ultra thick, red] (6,-2.8)--(6.9, -1.9) (7,-1.8)--(7.9, -2.7) (8,-2.8)--(8.9, -3.7); 

\draw [scale=0.8, ultra thick, red] (13,-3.8)--(13.9, -2.9) (14,-2.8)--(14.9, -3.7) (15,-3.8)--(15.9,-4.7) ;

\end{tikzpicture}
\end{center} 
\caption{A leaf-marked binary tree,  its corresponding full binary tree, and the Dyck path with marked 
$UDD$-patterns.} \label{Fig:tree-to-path} 
\end{figure}
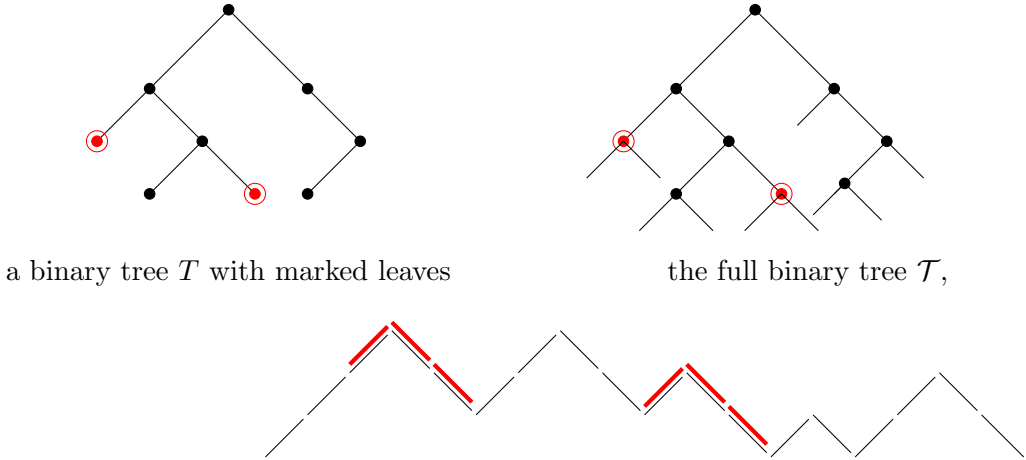 

Figure \ref{Fig:tree-to-path} shows a leaf-marked binary tree , its corresponding full binary tree 
and Dyck path with marked $UDD$-patterns. 

Let $a_{n,k}$ be the number of binary trees on $n$ vertices with $k+1$ leaves, and $P_{n,k}$ be 
the number of leaf-marked binary trees on $n$ vertices with $k$ marked leaves. 
Then the arrays $\{a_{n,k}: n \geq 1, 0\leq k \leq \lfloor \frac{n-1}{2} \rfloor \}$ and   
$\{P_{n,k}: n \geq 1, 0\leq k \leq \lfloor \frac{n-1}{2} \rfloor \}$ 
are closely related to the Catalan's and Borel's triangles.  In particular, $P_{n,k}=B_{n-1-k,k}$. 
The entries of $\{a_{n,k}\}$ and  $\{P_{n,k}\} $ for small values of $n$ and $k$ are listed below.  

\begin{equation*} 
\begin{array}{ccc} 
\begin{tabular}{c|cccc   } 
 $a_{n,k}$ & k=0  & 1 & 2     & 3   \\  \hline 
 n=0                        & 1  &    &        &      \\ 
 1                        & 1 &    &        &      \\
 2                        & 2  &  &        &      \\
 3                       & 4 & 1 &        &       \\
 4                       & 8  & 6 &      &       \\
 5                       & 16& 24 &2  &      \\
 6                      & 32 & 80 & 20 & 
 \end{tabular}

& \hspace{2cm}  &

\begin{tabular}{c|cccc   } 
 $P_{n, k}  $ & k= 0  & 1 & 2     & 3   \\  \hline 
 n=1                        & 1  &    &        &      \\ 
 2                        & 2  &    &        &      \\
 3                        & 5  & 1 &        &      \\
 4                       & 14 & 6 &        &       \\
 5                       & 42 & 28& 2    &       \\
 6                       & 132&120&20 &      \\
 7            &  429 & 495  & 135  & 5  
           \end{tabular}
           \end{array} 
\end{equation*}

% This part is done after the bijection

The array $\{a_{n,k}\}$ is  the sequence A091894 in OEIS, whose values are given by
\[
a_{n,k} = 2^{n-2k-1} \binom{n-1}{2k} C_k 
\]
for $n \geq 1$ and  whose row-sum is the Catalan sequence. The array $\{P_{n,k}\}$ can be obtained 
from Borel's triangle by reading the entries diagonally from lower left to upper right. 
Let $\mathcal{A}(t,x)=\sum_{n \geq 1, k\geq 0} a_{n,k} t^k x^n $ and 
$\mathcal{P}(t,x)= \sum_{n \geq 1, k \geq 0} P_{n,k} t^k x^n $ be the bivariate generating functions, then 
we have 
\[
\mathcal{P}(t,x) = \mathcal{A}(1+t, x) \quad  \text{and} \qquad \mathcal{P}(t,x)= x\mathcal{B}(tx, x). 
\] 
From these equations we can derive 
\[
\mathcal{A}(t,x) = x\mathcal{C} (1-x+tx,x),
\]
where $\mathcal{C}(t,x)$ and $\mathcal{B}(t,x)$ are the bivariate generating functions for Catalan's and 
Borel's triangles, respectively. 
Explicitly, one can compute that 
\[
\mathcal{A}(t,x) = \frac{1-2x-\sqrt{1-4x+4x^2-4x^2t} } {2tx} , 
\] 
which can also be  derived directly from the binary tree structures.

\subsection{Matchings Avoiding Pair of Patterns} \label{subsec: Pattern} 

%\input{3matching.tex}
%input file for pattern avoiding matchings

In~\cite[Section 6]{BE13}, Bloom and Elizalde gave nice results on pattern avoiding matchings. 
In particular, the number of pattern avoiding matchings from their class I is the generalized Catalan number. 
This suggests that there may be certain combinatorial statistics on those matchings that are enumerated by $B_{n,k}$. In this subsection we 
give explicit description of such statistics. 

For $2n$ points $1,2,\ldots, 2n$ on a horizontal line, we represent a \emph{matching} on these $2n$ points by drawing $n$ arcs between pairs of points.
A \emph{crossing} is formed by two arcs $(i,j)$ and $(k,\ell)$ if $i<k<j<\ell$, and a \emph{nesting} 
is formed by two arcs $(i,j)$ and $(k,\ell)$ where $i<k<\ell<j$. Denote by $\mathcal{M}_n$ the set of 
all matchings on $2n$ elements. For a matched pair $(i,j)$ with $i<j$, we call $i$ an \emph{opener} 
and $j$ a \emph{closer}.

A \emph{full rook placement} is a pair $(R,F)$ where $F$ is a Ferrers board and $R$ is a set of rooks 
 placed in $F$ such that each row and column of $F$ contains exactly one rook. 
Denote by $\mathcal{R}_F$ the set of all full rook placements on a given $F$, 
and $\mathcal{R}_n$ the set of all full rook placements on Ferrers boards having $n$ nonempty rows and 
columns. Throughout this section, a Ferrers board is left- and bottom-aligned, 
and we will use the name \emph{rook placement} to refer to a full rook placement.

\begin{figure}[htpb] 
\begin{center}
\begin{tikzpicture}[line cap=round,line join=round,>=triangle 45,x=1.0cm,y=1.0cm,scale=.7]
\clip(-.5,-1.5) rectangle (19,4.5);
%matching
\draw [fill] (0,0) circle [radius=.1];
\draw [fill] (1,0) circle [radius=.1];
\draw [fill] (2,0) circle [radius=.1];
\draw [fill] (3,0) circle [radius=.1];
\draw [fill] (4,0) circle [radius=.1];
\draw [fill] (5,0) circle [radius=.1];
\draw [fill] (6,0) circle [radius=.1];
\draw [fill] (7,0) circle [radius=.1];
\draw [fill] (8,0) circle [radius=.1];
\draw [fill] (9,0) circle [radius=.1];

\draw [line width=1.2pt] (0,0) to[bend left] (4,0);
\draw [line width=1.2pt] (1,0) to[bend left] (2,0);
\draw [line width=1.2pt] (3,0) to[bend left] (9,0);
\draw [line width=1.2pt] (5,0) to[bend left] (7,0);
\draw [line width=1.2pt] (6,0) to[bend left] (8,0);

\node at (0,-.5) {1};
\node at (1,-.5) {2};
\node at (2,-.5) {3};
\node at (3,-.5) {4};
\node at (4,-.5) {5};
\node at (5,-.5) {6};
\node at (6,-.5) {7};
\node at (7,-.5) {8};
\node at (8,-.5) {9};
\node at (9,-.5) {10};

\node at (11,1) {$\longrightarrow$};

%rook placement
\draw [line width=1.2pt] (13,4)--(13,-1)--(18,-1)--(18,2)--(16,2)--(16,3)--(15,3)--(15,4)--(13,4);

\draw [line width=1.2pt] (13,3)--(15,3) (13,2)--(16,2) (13,1)--(18,1) (13,0)--(18,0) (14,4)--(14,-1) (15,3)--(15,-1) (16,2)--(16,-1) (17,2)--(17,-1);
\draw [line width=1.2pt] (14,3)--(13,2) (13,3)--(14,2) (15,4)--(14,3) (14,4)--(15,3) (15,0)--(16,-1) (16,0)--(15,-1) (16,2)--(17,1) (17,2)--(16,1) (17,1)-- (18,0) (18,1)--(17,0);

\node at (13.5,4.3) {1};
\node at (14.5,4.3) {2};
\node at (15.3,3.5) {3};
\node at (15.7,3.3) {4};
\node at (16.3,2.5) {5};
\node at (16.7,2.3) {6};
\node at (17.5,2.3) {7};
\node at (18.3,1.5) {8};
\node at (18.3,.5) {9};
\node at (18.3,-.5) {10};
\end{tikzpicture}
\end{center}
\caption{A bijection between matchings and rook placements.}\label{Fig:matching-to-rook}
\end{figure}
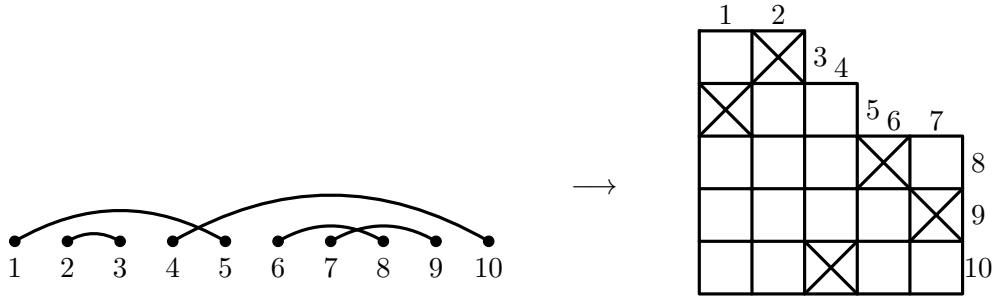

The sets $\match_n$ and $\mathcal{R}_n$ are related by the following natural bijection, see~\cite{deMier07, Jelinek07}. Given a rook placement $F\in \mathcal{R}_n$, label its border from top-left to lower-right by numbers $1, 2,\ldots, 2n$. If a rook has column labeled $s$ and row labeled $t$, draw an arc that connects points $s$ and $t$ in the matching. Conversely, for a matching $M\in \match_n$, draw a path from $(0,n)$ to $(n,0)$ by reading the vertices of $M$ in increasing order, and each opener corresponds to an east step while each closer corresponds to a south step. The resulting path will be the border of the Ferrers board $F$. A matched pair of $(s,t)$ then corresponds to a rook in $F$ that appears in column associated to $s$ and row associated to $t$. See the example in Figure~\ref{Fig:matching-to-rook} for an illustration.

To better describe the location of a rook, we consider the lower-left corner of a Ferrers board $F$ as the origin $(0,0)$ and label the columns of $F$ from left to right as $1,2, \ldots, n$. The rows of $F$ are labeled from bottom to top in the same manner. Then we say a rook is in the position $(i,j)$ if it appears in column $i$ and row $j$. Note that this labeling is different from the one in Figure~\ref{Fig:matching-to-rook}.

With above bijection and coordinates for rooks, we can translate the terms nesting and crossing to rook placements.

\begin{defi}
Given a rook placement $(R,F)$, if two rooks appear in positions $(i,j)$ and $(k,\ell)$, then they form a nesting pair if $i<k$ and $j <\ell$, that is, the rook $(i,j)$ is to the southwest of rook $(k,\ell)$. They form a crossing pair if $i>k$ and $j<\ell$, that is, the rook $(i,j)$ is to the southeast of rook $(k, \ell)$ and $F$ contains a square in the position $(i,\ell)$.
\end{defi}

Next we define pattern avoiding matchings and rook placements.

\begin{defi}
Let $M\in \match_n$ be a matching, then it avoids a pattern $\sigma\in \mathfrak{S}_k$ if there are no $2k$ vertices $i_1<i_2<\cdots<i_{2k}$ such that $M$ contains all pairs of $(i_r, i_{2k+1-\sigma(r)})$ for $1\leq r\leq k$.

Equivalently, a  rook placement  $(R, F)$ avoids a pattern $\sigma \in \mathfrak{S}_k$ if it does not contain a square subboard of size $k \times k$ such that the restriction of $F$ on the subboard is 
$\{(r, \sigma(r)): 1\leq r \leq k\}$. 
\end{defi}

%Pattern avoiding rook placements can be defined in a similar way and is omitted here.

For a pair of patterns $\{\sigma, \tau\}$  let $\match_n(\sigma, \tau)$ be the set of  all matchings of $[2n]$  that avoid the patterns $\sigma$ and $\tau$, and $\mathcal{R}_n(\sigma, \tau)$ be the corresponding set of rook placements.  
Then we have:

\begin{thm}
The entry $B_{n,k}$ of Borel's triangle enumerates the following sets.
\begin{enumerate}[(1)]
\item $(123, 213)$-avoiding matchings on $2(n+1)$ elements with $k$ nestings.
\item $(231, 321)$-avoiding matchings on $2(n+1)$ elements with $k$ crossings.
\end{enumerate}
\end{thm}

\begin{figure}[htpb]
\begin{center}
\begin{tikzpicture}[scale=.5]
\draw [line width=1.2pt] (0,0)-- (0,3)--(3,3)--(3,0)--(0,0);

\draw [line width=1.2pt] (1,3)--(1,0) (2,3)--(2,0) (3,3)--(3,0) (0,1)--(3,1) (0,2)--(3,2) (0,3)--(3,3);

\draw [line width=1.2pt] (0,0)--(1,1) (1,0)--(0,1) (1,1)--(2,2) (2,1)--(1,2) (2,2)--(3,3) (3,2)--(2,3);
\node at (1.5,-.5) {$123$};

\begin{scope}[shift={(6,0)}]
\draw [line width=1.2pt] (0,0)-- (0,3)--(3,3)--(3,0)--(0,0);

\draw [line width=1.2pt] (1,3)--(1,0) (2,3)--(2,0) (3,3)--(3,0) (0,1)--(3,1) (0,2)--(3,2) (0,3)--(3,3);

\draw [line width=1.2pt] (0,1)--(1,2) (0,2)--(1,1) (1,0)--(2,1) (1,1)--(2,0) (2,2)--(3,3) (3,2)--(2,3);
\node at (1.5,-.5) {$213$};
\end{scope}
\end{tikzpicture}
\end{center}
\caption{Rook placements of patterns $123$ and $213$ respectively.}
\label{Fig:rook_pattern}
\end{figure}
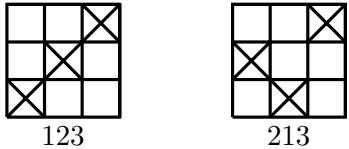

Figure~\ref{Fig:rook_pattern} gives an example of rook placements of patterns $123$ and $213$ respectively.

\begin{proof}
We prove (1) by constructing a bijection between rook placements in $\mathcal{R}_{n+1}(123,213)$ with $k$ nesting pairs  and Dyck paths of semilength $n+1$ 
with $k$ marked down-steps that are not at ground level.

For any rook placement $F\in \mathcal{R}_{n+1}(123, 213)$, its boundary is a Dyck path rotated $45^o$, 
so that down-steps are south steps on the boundary of $F$. To mark this Dyck path, if $r_1 = (i,j)$ and $r_2 = (k,\ell)$ is a nesting pair where $i<k$ and $j<\ell$,  mark the south step of the Dyck path corresponding to row $\ell$. Since the rook placement is $123$- and $213$-avoiding, each $r_2$ corresponds to exactly one $r_1$, that is, $r_2$ contributes to only one nesting pair and a down-step will be marked at most once.
Thus $k$ nesting pairs of rooks give $k$ marked down-steps in the Dyck path.

Next we need to show these marks are not at ground level. Assume the mark of row $\ell$ corresponds to a down-step at ground level, that is, this step has the lower endpoint on the diagonal $x+y=n$. Consider the subboard $F'$ obtained by deleting all rows below $r_2$. Since $r_2$ lies in a row on the diagonal, $F'$ contains $m$ columns and $m$ rows. Since $F$ is a full rook placement, $F'$ must be a full rook placement. But $r_1$ is to the southwest of $r_2$, this means at least one column in $F'$ must be empty, a contradiction. Thus $r_2$ can not be a rook in a row on the diagonal, and the resulting Dyck path is one where all marked down-steps are not at ground level.

\begin{figure}[htpb]
\begin{center}
\begin{tikzpicture}[line cap=round,line join=round,>=triangle 45,x=1.0cm,y=1.0cm,scale=.4]
\clip(-.5,-7.5) rectangle (26,10.5);
%nesting rook placement
\draw [line width=1.2pt] (0,0)-- (0,10)--(4,10)--(4,7)--(7,7)--(7,6)--(9,6)--(9,2)--(10,2)--(10,0)--(0,0);

\draw [line width=1.2pt] (1,10)--(1,0) (2,10)--(2,0) (3,10)--(3,0) (4,7)--(4,0) (5,7)--(5,0) (6,7)--(6,0) (7,6)--(7,0) (8,6)--(8,0) (9,2)--(9,0) (0,1)--(10,1) (0,2)--(9,2) (0,3)--(9,3) (0,4)--(9,4) (0,5)--(9,5) (0,6)--(7,6) (0,7)--(4,7) (0,8)--(4,8) (0,9)--(4,9);

%rooks
\draw [line width=1.2pt] (1,10)--(2,9) (2,10)--(1,9) (2,9)--(3,8) (3,9)--(2,8) (3,8)--(4,7) (4,8)--(3,7) (4,7)--(5,6) (5,7)--(4,6) (0,6)--(1,5) (1,6)--(0,5) (5,3)--(6,2) (6,3)--(5,2) (6,5)-- (7,4) (7,5)--(6,4) (7,4)--(8,3) (8,4)--(7,3) (8,2)--(9,1) (9,2)--(8,1) (9,1)--(10,0) (10,1)--(9,0);

\node at (7.5,-1) {a $(123,213)$-avoiding rook placement};

%marked Dyck path
\draw [line width=1.2pt] (15,0)-- (15,10)--(19,10)--(19,7)--(22,7)--(22,6)--(24,6)--(24,2)--(25,2)--(25,0)--(15,0);

\draw [line width=1.2pt] (16,10)--(16,0) (17,10)--(17,0) (18,10)--(18,0) (19,7)--(19,0) (20,7)--(20,0) (21,7)--(21,0) (22,6)--(22,0) (23,6)--(23,0) (24,2)--(24,0) (15,1)--(25,1) (15,2)--(24,2) (15,3)--(24,3) (15,4)--(24,4) (15,5)--(24,5) (15,6)--(22,6) (15,7)--(19,7) (15,8)--(19,8) (15,9)--(19,9);

\draw [line width=1.2pt, color=red] (16,10)--(17,9) (17,10)--(16,9) (17,9)--(18,8) (18,9)--(17,8) (18,8)--(19,7) (19,8)--(18,7) (19,7)--(20,6) (20,7)--(19,6) (21,5)-- (22,4) (22,5)--(21,4) (22,4)--(23,3) (23,4)--(22,3);

\draw [line width=1.2pt] (15,6)--(16,5) (16,6)--(15,5) (20,3)--(21,2) (21,3)--(20,2) (23,2)--(24,1) (24,2)--(23,1) (24,1)--(25,0) (25,1)--(24,0);

%marks
\draw [line width=1.2pt,color=red] (19.3,9.9)--(19.3,9.1) (19.3,8.9)--(19.3,8.1) (19.3,7.9)--(19.3,7.1) (22.3,6.9)--(22.3,6.1) (24.3,4.9)--(24.3,4.1) (24.3,3.9)--(24.3,3.1);

\node at (20,-1) {nestings and marks};

%Dyck path
\draw [line width=1.2pt] (4,-7)--(4.9,-6.1) (5,-6)--(5.9,-5.1) (6, -5)--(6.9, -4.1) (7, -4)--(7.9, -3.1) (8,-3)--(8.9,-3.9) (9, -4)--(9.9,-4.9) (10,-5)--(10.9,-5.9) (11,-6)--(11.9,-5.1) (12,-5)--(12.9,-4.1) (13,-4)--(13.9,-3.1) (14,-3)--(14.9,-3.9) (15.,-4.)--(15.9,-3.1) (16,-3)--(16.9,-2.1) (17,-2)--(17.9,-2.9) (18,-3)--(18.9,-3.9) (19,-4)--(19.9,-4.9) (20,-5)--(20.9,-5.9) (21,-6)--(21.9,-5.1) (22,-5)--(22.9,-5.9) (23,-6)--(24,-7);

\draw [line width=1.2pt, color=red] (8.2,-2.8)--(9.1,-3.7) (9.2,-3.8)--(10.1, -4.7) (10.2, -4.8)--(11.1,-5.7) (14.2, -2.8)--(15.1,-3.7) (18.2,-2.8)--(19.1, -3.7) (19.2,-3.8)--(20.1,-4.7);

%\node at (13, -8) {marked Dyck path};
\end{tikzpicture}
\end{center}
\caption{An example of $(123, 213)$-avoiding rook placement and its corresponding Dyck path where the marked down-steps are not at ground level.}\label{Fig:nesting-to-path}
\end{figure}

To see the inverse, consider a Dyck path as the border of a Ferrers board $F$. Place rooks in $F$ from the top row to bottom in the following manner: If a row corresponds to an unmarked down-step in the Dyck path, place a rook in the leftmost empty column. If a row corresponds to a marked down-step, place a rook in the second leftmost empty column in $F$ so that there is one empty column to its left. Since this down-step is not at ground level, there are more up-steps than down-steps in the segment of the Dyck path before this step, that is, there are more columns than rows in the Ferrers board prior to this row. Since we place rooks by rows, there is always an empty column to the left of a rook corresponding to a marked down-step. 

Since the rooks are placed from top to bottom, for any rook $r$ in this placement, there is at most one rook that is to the southwest of it, hence it is $(123,213)$-avoiding. Moreover, the rooks $r_1 = (i,j)$ and $r_2 = (k,\ell)$ with $i<k$ is a nesting pair if and only if $r_2$ is in a marked row and $r_1$ is in the first unmarked row below it, thus marks in the Dyck path correspond to nesting pairs in the rook placement and hence this is a bijection. See Figure~\ref{Fig:nesting-to-path} for an example.

The proof of (2) is similar, we find a bijection between $(231, 321)$-rook placements and Dyck path with marked up-steps that are not at ground level.

For any rook placement $F\in \mathcal{R}_{n+1}(231, 321)$, its border is a rotated Dyck path in which up-steps are east-steps. If two rooks $r_1 = (i,j)$ and $r_2 = (k,\ell)$ is a crossing pair and $r_1$ is to the northwest of $r_2$, mark the east step of the Dyck path corresponding to the column of $r_2$. Conversely, given a Dyck path with marked up-steps, use it as the border of a Ferrers board $F$. Place the rooks in $F$ from right to left in the following manner: If the column corresponding to an east step is not marked, place a rook in the highest empty row. If the column corresponding to an east step is marked, place a rook in the second highest empty row such that there is an empty row above the rook.

%Assume this east step is on the ground, then consider the Ferrers board $F'$ obtained by all columns to the right of $r_2$ and all rows except the one $r_2$ lies in. $F'$ contains $m$ rows and $m$ columns, but since there is always a square in the position $(k, j)$ by definition of crossing, only $m-1$ rooks are allowed to be placed in $F'$ by rows while $m$ rooks need to be placed in $F'$ by columns, as $r_1$ is placed in the row corresponding to $(k,j)$, a contradiction. Thus this step in the Dyck path is not at ground level.

\begin{figure}[htpb]
\begin{center}
\begin{tikzpicture}[line cap=round,line join=round,>=triangle 45,x=1.0cm,y=1.0cm,scale=.4]
\clip(-.5,-.5) rectangle (26,17.5);
%Dyck path
\draw [line width=1.2pt] (4,12)--(4.9,12.9) (5,13)--(5.9,13.9) (6,14)--(6.9, 13.1) (7,13)--(7.9,13.9) (8,14)--(8.9,14.9) (9,15)--(9.9,15.9) (10,16)--(10.9,16.9) (11,17)--(11.9,16.1) (12,16)--(12.9,15.1) (13,15)--(13.9,15.9) (14,16)--(14.9,15.1) (15,15)--(15.9,14.1) (16,14)--(16.9,13.1) (17,13)--(17.9,13.9) (18,14)--(18.9,14.9) (19,15)--(19.9,15.9) (20,16)--(20.9,15.1) (21,15)--(21.9,14.1) (22,14)--(22.9,13.1) (23,13)--(23.9,12.1);

\draw [line width=1.2pt,color=red] (7.9,14.2)--(8.8,15.1) (8.9,15.2)--(9.8,16.1) (12.9,15.2)--(13.8,16.1) (16.9,13.2)--(17.8,14.1) (17.9,14.2)--(18.8,15.1) (18.9,15.2)--(19.8,16.1);

\node at (13,11) {a Dyck path with marked non-ground level up-steps};

%marked border
\draw [line width=1.2pt] (0,0)--(10,0)--(10,4)--(7,4)--(7,7)--(6,7)--(6,9)--(2,9)--(2,10)--(0,10)--(0,0);
\draw [line width=1.2pt] (0,1)--(10,1) (0,2)--(10,2) (0,3)--(10,3) (0,4)--(7,4) (0,5)--(7,5) (0,6)--(7,6) (0,7)--(6,7) (0,8)--(6,8) (0,9)--(2,9) (1,10)--(1,0) (2,9)--(2,0) (3,9)--(3,0) (4,9)--(4,0) (5,9)--(5,0) (6,7)--(6,0) (7,4)--(7,0) (8,4)--(8,0) (9,4)--(9,0);

%marks
\draw [line width=1.2pt,color=red] (3.1,9.3)--(3.9,9.3) (4.1,9.3)--(4.9,9.3) (6.1,7.3)--(6.9,7.3) (7.1,4.3)--(7.9,4.3) (8.1,4.3)--(8.9,4.3) (9.1,4.3)--(9.9,4.3);

%rooks
\draw [line width=1.2pt,color=red] (9,3)--(10,2) (10,3)--(9,2) (8,2)--(9,1) (9,2)--(8,1) (7,1)--(8,0) (8,1)--(7,0) (6,6)--(7,5) (7,6)--(6,5) (4,7)--(5,6) (5,7)--(4,6) (3,5)--(4,4) (4,5)--(3,4);

\draw [line width=1.2pt] (5,9)--(6,8) (6,9)--(5,8) (2,8)--(3,7) (3,8)--(2,7) (1,10)--(2,9) (2,10)--(1,9) (0,4)--(1,3) (1,4)--(0,3);

%(231,321)-avoiding rook placement
\draw [line width=1.2pt] (15,0)--(25,0)--(25,4)--(22,4)--(22,7)--(21,7)--(21,9)--(17,9)--(17,10)--(15,10)--(15,0);
\draw [line width=1.2pt] (15,1)--(25,1) (15,2)--(25,2) (15,3)--(25,3) (15,4)--(22,4) (15,5)--(22,5) (15,6)--(22,6) (15,7)--(21,7) (15,8)--(21,8) (15,9)--(17,9) (16,10)--(16,0) (17,9)--(17,0) (18,9)--(18,0) (19,9)--(19,0) (20,9)--(20,0) (21,7)--(21,0) (22,4)--(22,0) (23,4)--(23,0) (24,4)--(24,0);

%rooks
\draw [line width=1.2pt] (24,3)--(25,2) (25,3)--(24,2) (23,2)--(24,1) (24,2)--(23,1) (22,1)--(23,0) (23,1)--(22,0) (21,6)--(22,5) (22,6)--(21,5) (19,7)--(20,6) (20,7)--(19,6) (18,5)--(19,4) (19,5)--(18,4);

\draw [line width=1.2pt] (20,9)--(21,8) (21,9)--(20,8) (17,8)--(18,7) (18,8)--(17,7) (16,10)--(17,9) (17,10)--(16,9) (15,4)--(16,3) (16,4)--(15,3);

\end{tikzpicture}
\end{center}
\caption{An example of $(231,321)$-avoiding rook placement and its corresponding Dyck path with marked up-steps that are not at  ground level.}\label{Fig:crossing-to-path}
\end{figure}
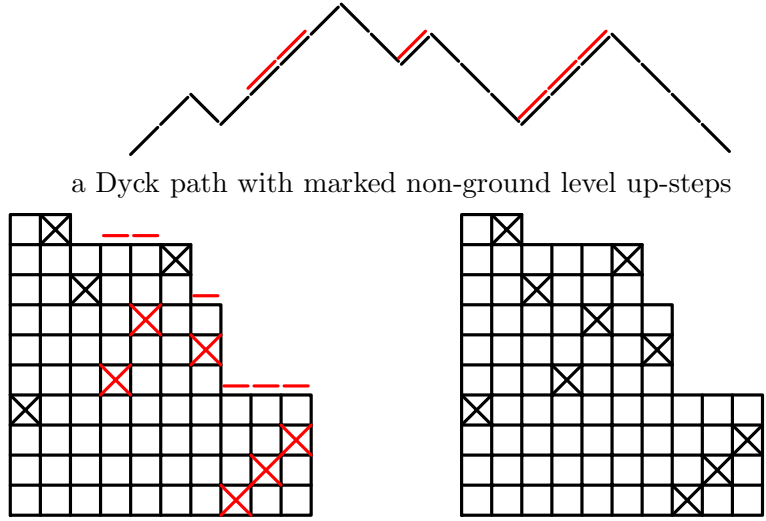

%Since the labeled east step is not at ground level, there are more down-steps than up-steps in the segment after this step, that is, there are more rows than columns in the Ferrers board, since we place rooks by columns, there must exist an empty row to the north of such a rook.
By a similar argument as in the proof of (1), we can show that this map is well-defined and bijective, moreover, the resulting rook placement is $(231,321)$-avoiding. The details are left to the readers. See Figure~\ref{Fig:crossing-to-path} for an example.
\end{proof}

In~\cite[Section 6]{BE13} the authors gave seven other families of pattern avoiding matchings that 
are enumerated by the generalized Catalan number. 
For each family, we can define a statistic $s(M)$ for $M \in \mathcal{M}_{n+1}(\sigma, \tau)$ such that the subset $S_{n,m}:=\{ M\in \mathcal{M}_{n+1}(\sigma, \tau): s(M)=m\}$ has cardinality $B_{n,m}$.
We describe in each case the statistic $s(M)$ in Table~\ref{Table:patterns} without proofs. 

\begin{table}[htpb]
\begin{center}
\begin{tabular}{|l|l|}
\hline
Patterns  $(\sigma, \tau)$ & $s(M)$ for $M\in \mathcal{M}_{n+1}(\sigma, \tau)$\\
\hline
\hline
$\{231, 312\}$ & $|\{(i,j): \text{there exists some pair } (k,\ell) \text{ such that } k<i<j<\ell, $\\
					& $\text{ and there is no pair } (r,s) \text{ where } r<i<s<j.\}|$\\
\hline
$\{132,213\}$ & $|\{(i,j): \text{there exists some pair } (k,\ell) \text{ such that } k<i<j<\ell.\}|$\\
\hline
$\{213,231\}$ & $|\{(i,j): \text{there exists some pair } (k,\ell) \text{ such that } k<i<j<\ell,$\\
					& $\text{ and there is pair } (i-1, r) \text{ where } i<r<j.\}|$\\
\hline
$\{132, 312\}$ & $|\{(i,j): \text{there exists some pair } (k,\ell) \text{ such that } k<i<j<\ell.\}|$\\
\hline
$\{321,312\}$ & $|\{(i,j): \text{there exists some pair } (k,\ell) \text{ such that } i<k<j<\ell.\}|$\\
\hline
$\{213,312\}$ & $|\{(i,j): \text{there exists some pair } (k,\ell) \text{ such that } k<i<j<\ell,$\\
					& $\text{ and there is pair } (r, j+1) \text{ where } i<r<j.\}|$\\
\hline
$\{132,231\}$ & $|\{(i,j): \text{there exists some pair } (k,\ell) \text{ such that } k<i<j<\ell.\}|$\\
\hline
\end{tabular}
\caption{Seven cases of pattern avoiding matchings on $2(n+1)$ elements such that $|\{M\in \mathcal{M}_{n+1}(\sigma, \tau): s(M)=m\}| = B_{n,m}$.}\label{Table:patterns}
\end{center}
\end{table}

\section{Restricted Dumont Permutation and Marked Dyck paths} \label{Sec:Dumont}

%\input{4dumont.tex} 
% This is the input file of Section 4 on restricted Dumond permutations 

% Definition.  Results/ equation,  and bijection to Dyck paths. 
Dumont permutations were introduced by Dumont \cite{Dumont74}  as a combinatorial interpretation for Genocchi numbers, which is a multiple of  Bernoulli numbers.  Surprisingly, 
the sequence of generalized Catalan numbers also   appears in  the enumeration of 
restricted Dumont permutations.  Using generating functions 
Burnstein proved  in  \cite[Theorem 3.1]{Burnstein05} that the $n$-th generalized Catalan number  
counts the number of Dumont permutations of the first kind with  length $2n$ that 
avoid patterns $(2413, 3142)$. 
In this section we construct a bijection between those restricted Dumont permutations and Dyck paths with marked down-steps not at ground level.  
%by  analyzing the recurrence structures of the involved objects, 
Our construction  leads to a functional equation for the generalized Catalan numbers and Borel's triangle, 
which provides a useful tool to analyze the structures counted by Borel's triangle, (c.f.  Section \ref{Sec:functional} ).

\begin{defi} \label{Df:dumont} 
A \emph{Dumont permutation of the first kind} is a permutation $\pi$ of even length, i.e., 
$\pi \in \mathfrak{S}_{2n}$,  of which each even entry is followed by a descent and each odd entry is followed by an ascent or is at the end of the permutation. In other words, for every $i=1, 2, \dots, 2n$, 
\begin{eqnarray*} 
\pi(i) \text{ is even }  & \Longrightarrow  &  i < 2n \text{ and } \pi(i) > \pi(i+1),   \\
\pi(i) \text{ is odd }    & \Longrightarrow  &  \pi(i) < \pi(i+1) \text{ or } i=2n.
\end{eqnarray*} 
\end{defi} 

Denote by $\mathcal{D}_{2n}^1$ the set of Dumont permutations of the first kind of length $2n$. For example, $\mathcal{D}_2^1=\{ 21\}$ and 
$\mathcal{D}_4^1=\{ 2143, 3421, 4213 \}$.  Note that in a Dumont permutation of the first kind, 
$21$ must appear together, and the permutation always ends with an odd number. 

\textsc{Remark}. Dumont permutations of odd lengths can be defined similarly by the conditions in 
Definition \ref{Df:dumont}. To meet the condition on the odd entries, the largest entry $2n+1$ must 
appear at the end. Hence $\mathcal{D}^1_{2n+1}$ can be obtained  simply by adjoining  $2n+1$ to the end 
of each permutation in $\mathcal{D}^1_{2n}$ and $|\mathcal{D}^1_{2n+1} | = |\mathcal{D}^1_{2n}|$.

Let $S$ be a set of permutations and $ \mathcal{D}^1_{2n}(S )$ be the set of Dumont permutations of length $2n$ that avoid patterns in $S$, that is, there is no subsequence  that is order isomorphic to any pattern in $S$. 
Burnstein analyzed the structures and enumerated Dumont permutations avoiding  a certain pattern or sets of patterns of length 3 and 4. One of his results relates to the $n$-th generalized Catalan number
$C(2,n)$ defined via equations.\eqref{eq2.2-2} and \eqref{eq2.2-3}.

\begin{thm}[Burnstein]  \label{thm:Burnstein} 
\[
|\mathcal{D}_{2n}^1 (2413, 3142) | = C(2,n). 
\]
\end{thm} 

For example, $\mathcal{D}^1_{4}(2413, 3142) =\mathcal{D}^1_4$, but $|\mathcal{D}^1_6|= 17$ while 
$|\mathcal{D}^1_{6}(2413, 3142)|=13$ since there are  four Dumont permutations of the first kind and length 6 that contain patterns $2413$ or $3142$. More precisely, they  are $\{ \mathbf{36}42\mathbf{15}, 
\mathbf{4}2\mathbf{15} 6\mathbf{3},  \mathbf{4}3\mathbf{62}1\mathbf{5},  \mathbf{4}21\mathbf{635} \}$,
where an occurrence of the  forbidden patterns is boldfaced.

%Burstein proved Theorem \ref{thm:Burnstein} using funcitonal equations. 
We will prove Theorem \ref{thm:Burnstein} by constructing a bijection $\rho$  from the set 
$\mathcal{D}^1_{2n}(2413, 3142)$ to the set of Dyck paths of semi-length $n$ with marked 
down-steps not at ground level. 
First we need to understand the structure of the permutations in  $\mathcal{D}^1_{2n}(2413, 3142)$. 
The following decomposition is also the idea behind Burnstein's generating function argument.  
To construct the desired bijection, we prove more properties for the decomposition and correct a 
typo\footnote{The block form \eqref{dumont-form} appeared in \cite{Burnstein05}, but the term $(2k)$ was written 
as $(2k-2)$.} 
in \cite{Burnstein05}. 

Let $\pi \in \mathcal{D}^1_{2n}(2413, 3142)$.  Assume the last entry of $\pi$ is  $2k-1$ for some $k \in [1,n]$. Then $\pi$ can be decomposed into blocks of the form
\begin{equation} \label{dumont-form} 
 \cdots  A_i  B_i A_{i-1} B_{i-1} \cdots A_1 B_0  (2k) A_0  (2k-1), 
\end{equation}
where  $A_0$ is the part between entries $2k$ and $2k-1$, and $A_i$, $B_i$ are blocks of consecutive entries such that each entry in $A_i \  (i \geq 1) $ is less  than $2k-1$ and each entry in $B_i \ (i \geq 0)$ is 
greater than $2k$.  The blocks $A_1, B_1 \cdots $ are nonempty but $A_0$ and $B_0$ may be empty. 

\begin{lemma} 
 If $A_0$ is nonempty, then each entry in $A_0$ is less than $2k-1$. 
\end{lemma} 
\begin{proof} 
If $A_0$ is nonempty, since $2k$ is followed by a descent, the entry $i$  immediately after $2k$ is less than $2k-1$. If there is another  entry $j$  in $A_0$ that is larger than $2k$, then the four entries 
$2k, i , j, 2k-1$ form a pattern of 3142, a contradiction. 
\end{proof} 

\begin{lemma} 
\begin{enumerate} 
\item For any $i \geq 0$, assume  $ a  \in A_i$ and $b\in A_{i+1}$.   Then $a > b$. 
\item  For any $i \geq 0$, assume  $a \in B_i $ and $b \in B_{i+1}$.Then 
$a < b$. 
\end{enumerate} 
\end{lemma} 
\begin{proof} 
For (1), between $a$ and $b$ there is always an element $c $ which is larger than or equal to $2k$. 
If $b  >a$, then the terms $a, c, b$ and the last entry $2k-1$ of $\pi$ form a pattern of $2413$, a 
contradiction. 
Part (2) is similar, 
\end{proof} 

\begin{lemma} \label{lem:3} 
All the sets $|A_i|$ and $|B_j|$ are of even size. 
\end{lemma} 
\begin{proof} From the previous lemma each of $A_i$ and $B_j$ contains a set of consecutive integers. Lemma \ref{lem:3} follows from the fact that in each $A_i$, the smallest element is followed by an ascent, hence it must be odd. Similarly in each $B_i$, the largest element is followed by a descent  and hence must be even. 
\end{proof} 

For a string $W$ of distinct integers,  let $\tau(W)$ be the permutation of length $|W|$ that is order isomorphic to $W$.  For a permutation $\pi=a_1a_2 \cdots a_n$, let $c(\pi)$ be the complement of $\pi$, i.e., the permutation $b_1b_2\cdots b_n$ where $b_i=n+1-a_i$.

\begin{lemma}
 For all  $i$, $\tau(A_i)$  and $c(\tau(B_i))$ are 
  Dumont permutations of the first kind  that avoid $\{2413, 3142\}$. 
\end{lemma} 
\begin{proof} 
The statement for $A_i$ is obvious. For each $B_i$,  the sequence avoids $\{2413, 3142\}$, has an even length, and ends with an even entry. Hence taking the complement of $\tau(B_i)$, we obtain 
a Dumont permutation of the first kind that avoids $\{c(2413), c(3142)\} = \{3142, 2413\}$. 
\end{proof} 

Conversely, one checks that any permutation of the form \eqref{dumont-form} and satisfies the above lemmas is a Dumont permutation of the first kind avoiding $\{2413, 3142\}$. 

Now we are ready to describe the bijection $\rho$  from $\mathcal{D}^1_{2n}(2413, 3142)$ to the set of Dyck paths of semi-length $n$ with marked down-steps not at ground level. Our construction is defined recursively.

For $n=1$, 
the only permutation in $\mathcal{D}^1_2$ is $21$. Let  $\rho(21)$  be the unique Dyck path $UD$ of semi-length 1.

For $n \geq 2$, assume that $\pi \in \mathcal{D}^1_{2n}(2413, 3142)$ has form 
\eqref{dumont-form}. 
\begin{itemize}
\item \textbf{[Case 1]} If $\pi$ is of the form $ (2n) A_0 (2n-1)$, then $\rho(\pi)= U \rho(A_0) D$. 
\item \textbf{[Case 2]}  If the first non-empty block (from left) of $\pi$ is $B_i$, then $\rho(\pi)$ is the juxtaposition of 
$\rho(c(\tau(B_i)))$ and $\rho(\pi-B_i)$, where $\pi-B_i$ is the permutation obtained from $\pi$ by removing $B_i$. 
\item \textbf{[Case 3]} If the first non-empty block of $\pi$ is $A_i $ for some $i \geq 1$, 
then  let 
\begin{equation} \label{eq:sqcup} 
\rho(\pi) = \rho(\tau(A_i)) \sqcup \rho( \tau( \pi-A_i))
\end{equation} 
where $\pi-A_i$ is the permutation obtained from $\pi$ by removing $A_i$, and for two Dyck paths 
$P_1$ and $P_2$ with marked steps, the operation $P_1 \sqcup P_2$ is a new Dyck path with marked steps defined by the following rules. 
\begin{enumerate}[(i)]
\item Add a mark to the last down-step of $P_1$ to get $\tilde{P}_1$. 
\item Assume that   the primitive decomposition of path $P_2$ is $P_2=Q_1Q_2 \cdots Q_j$, where each $Q_i \ (1\leq i \leq j)$ is a primitive Dyck path, i.e., Dyck path that does not touch the $x$-axis except 
at the endpoints.  If $j \geq 2$, change $Q_{j-1} $ to $\tilde{Q}_{j-1}$ by add a mark to  the last down-step of $Q_{j-1}$. 
\item Assume that $Q_j=U R_j D$ where $Q_j$ is the last primitive factor of $P_2$. Set 
\[
P_1 \sqcup P_2 = U \tilde{P}_1 Q_1 \cdots Q_{j-2} \tilde{Q}_{j-1} R_j  D. 
\]
\end{enumerate}
\end{itemize} 

\begin{figure} [htpb] 
\begin{center} 
\begin{tikzpicture}[scale=.8]
\draw[thick] (0,1)--(1,2) (1.1,2.1)--(2,3)--(3,2) (3.1,1.9)--(4,1); 
\node at (2,0) {$\rho(4213)$}; 
\draw[thick] (5,1)--(6,2)--(7,1)--(8,2)--(9,1); 
\node at (7,0) {$\rho(3421)$}; 
\draw[thick] (10,1)--(11,2) (11.1,2.1)--(12,3)--(13,2) (13.1, 1.9)--(14,1); 
\draw[ultra thick,red] (12.2,3)--(13.2,2); 
\node at (12,0) {$\rho(2143)$}; 
\end{tikzpicture}
\caption{Bijection from $\mathcal{D}^1_4(2413, 3142)$ to marked Dyck paths. } \label{rho:n=4} 
\end{center}
\end{figure}
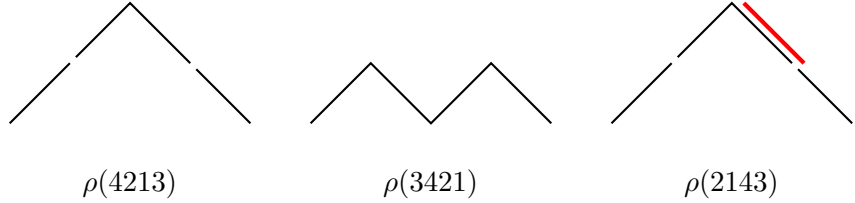

\begin{exa} 
For $n=2$, there are three permutations in $\mathcal{D}^1_{4}(2413,3142)$. The permutation
$4213$ is of Case 1 and $\rho(4213) = UUDD$. The permutation   $3421$ is of Case 2 and $\rho(3421)=UDUD$. The last permutation $2143$ is of Case 3, and hence 
$\rho(2143) = (UD) \sqcup (UD) = UU\hat{D}D$, where $\hat{D}$ means the corresponding down step is marked. See Figure  \ref{rho:n=4}, where the marked  down step is represented by a red line along the step. 
\end{exa} 

To see that $\rho$ is invertible and hence a bijection, we note some special properties of the preceding construction. 
\begin{prop} \label{prop:reverse} 
\begin{enumerate}
\item if $\pi$  is of the form $(2n) A_0 (2n-1)$, then $\rho(\pi)$ is a primitive Dyck path with no 
marked down-steps at level 0 or 1. 
\item If the first non-empty block of $\pi$ is $B_i$, then $\rho(\pi)$ is not primitive. 
\item If the first non-empty block of $\pi$ is $A_i$ for some $i \geq 1$, then $\rho(\pi)$ is a primitive Dyck path with at least one marked down-step at level 1. 
\begin{description}
\item  -- If $\rho(\pi)$ has  exactly one marked down-step at level 1, then $\pi$ is of the form $A_1 (2n) A_0 (2n-1)$.
\item -- If $\rho(\pi)$ has  at least two marked down-steps at level 1, then $B_{i-1} $ is not empty.  
\end{description}
\end{enumerate}
\end{prop} 

Given a Dyck path $P$ of length $2n$ with marked down-steps not at ground level, we can decide which case it comes from by using  Proposition~\ref{prop:reverse}. Then  the original permutation can be recovered by the recursive structure.  Case 1 is straightforward. For Case 2, the last primitive factor of $P$ corresponds to $\pi-B_i$, and the other part corresponds to $\tau(B_i)$, from which we can easily get $B_i$. For Case 3, the path $P$ is primitive. The part to the first marked down-step at level 1  allows us to recover  $A_i$, and the following part until the second marked down-step at level 1 (if there is one) allows us to  recover $P_2$
and hence  $\tau(\pi-A_i)$.  This proves that $\rho$ is invertible. 

We do not have a refinement of Theorem \ref{thm:Burnstein} in terms of Borel's triangle since we do not have  a clean description of the statistic on Dumont permutations 
that corresponds to the marked down-steps in Dyck paths.  
Nevertheless, the bijection $\rho$ allows us to count the number of primitive Dumont permutations. 

\begin{defi} 
 We say a Dumont permutation $\pi=\pi(1) \pi(2) \cdots \pi(2n)$  is \emph{decomposable} if there is an even number $2i$ with $i < n$ such that any term in $\pi(1), \pi(2), \cdots \pi(2i) $ is greater than any term in $\pi(2i+1), \dots, \pi(2n)$. In other words, $\pi(1), \dots, \pi(2i)$ is a permutation of  
 $2n-2i+1, \cdots, 2n$.   If a Dumont permutation is not decomposable, we say that it is \emph{primitive}. 
\end{defi} 

\begin{cor} 
The number of primitive Dumont permutations in $\mathcal{D}^1_{2n} (2413, 3142)$ is $2^{n-1} C_{n-1}$, where $C_{n-1}$ is the $(n-1)$-st Catalan number. 
\end{cor} 
\begin{proof} 
By Proposition~\ref{prop:reverse},  a permutation  $\pi$ in $\mathcal{D}^1_{2n} (2413, 3142)$ is primitive if and only if the image $\rho(\pi)$ is a primitive Dyck path with possibly some marked down-steps not at ground level. There are $C_{n-1}$ primitive Dyck paths, each with exactly $n-1$ down-steps not at ground level, hence $2^{n-1}$ ways to mark  those down-steps. 
\end{proof}

\section{Functional Equations of Borel's Triangle} \label{Sec:functional} 
%input file for tree decomposition

The bijection proof of Theorem \ref{thm:Burnstein} leads to a functional equation for the biviarate generating function of the entries in Borel's triangle. Let $f_{0,0}=1$ and $f_{n,k}$ be the number of Dyck paths of semi-length $n$ with $k$ marked down-steps not at ground level. We have  $f_{n,k}= B_{n-1,k}$ for $n \geq 1$  and 
$F(t,x)= 1+xB(t,x)$,   
where $F(t,x):=\sum_{n, k \geq 0} f_{n,k} t^k x^n  $ and  $B(t,x)$ is given in Eq.\eqref{eq-Borel2}. 

Using the block-decomposition \eqref{dumont-form} and writing $F(t,x)$ as $F$, we get
\begin{eqnarray}\label{eq-bivariate1} 
F= 1+ xF^2 + xF \sum_{i=1}^\infty t^{2i-1} \left( (F-1)^{2i-1} +(F-1)^{2i} \right)  \nonumber \\
+ xF\cdot t (F-1) \sum_{i=1}^\infty t^{2i-1} \left( (F-1)^{2i-1}+ (F-1)^{2i}\right).
\end{eqnarray}
In the above formula, $xF^2$ is the contribution of those Dumont permutations in 
$\mathcal{D}^1_{2n}(2413, 3142)$ which are of the form $B_0 (2k) A_0 (2k-1)$, where $B_0$ and $A_0$ may be empty. The first summation on the right-hand side is the contribution from all those  Dumont permutations with $B_0=\emptyset$ and $A_1 \neq \emptyset$, and the last summation is the contribution from those with $B_0 \neq \emptyset$ and $A_1 \neq \emptyset$. 
Simplifying Eq.\eqref{eq-bivariate1} we obtain the following proposition, 
where the special case when $t=1$ is given in \cite{Burnstein05}. 
\begin{thm} 
The bivariate generating function $F=F(t,x)$ satisfies the equation 
\begin{equation} \label{functional-eq} 
F= 1+xF^2 \cdot \frac{1}{1-t(F-1)}
\end{equation} 
\end{thm} 
It is easy to derive a quadratic equation from \eqref{functional-eq} 
\[
 (x+t)F^2 -(2t+1)F + t+1 =0,
\]
whose root with $F(0,0)=1$ is 
\begin{eqnarray} 
F(t,x) = \frac{ 2t+1-\sqrt{1-4x(1+t)}}{2(x+t)}  = \frac{1}{1-xC((1+t)x)}, 
\end{eqnarray}
which agrees with  \eqref{eq-Borel2}  by a simple algebraic manipulation.

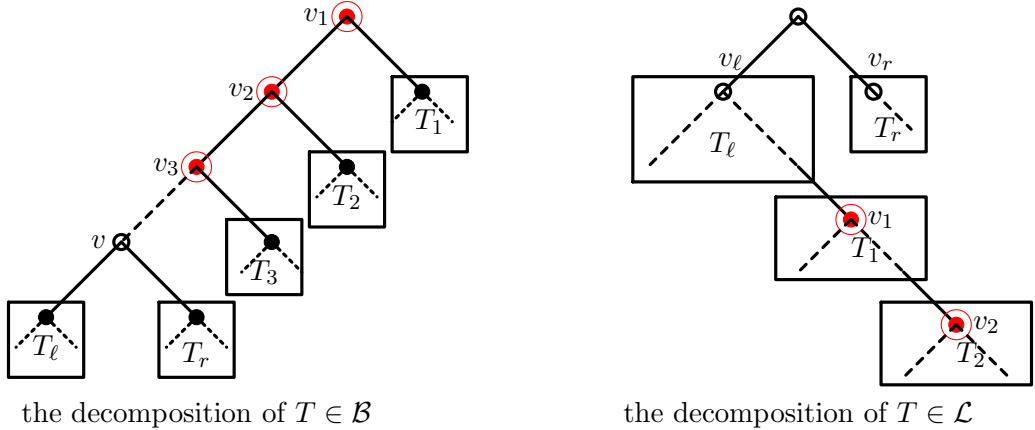
\begin{figure}[htbp] 
\begin{center}
\begin{tikzpicture}[line cap=round,line join=round,>=triangle 45,x=1.0cm,y=1.0cm]
\clip(2,-1) rectangle (17,5.5);

\begin{scope}[shift={(8,0)}]
%Not on right spine
\draw [line width=1.2pt] (5,5) circle [radius=.1];
\draw [line width=1.2pt] (5,5)--(6,4);
\draw [line width=1.2pt] (6,4) circle [radius=.1];
\draw [line width=1.2pt,dashed] (6,4)--(6.5,3.5);
\draw [line width=1.2pt] (5.7,4.2)--(6.7,4.2)--(6.7,3.2)--(5.7,3.2)--(5.7,4.2);
\node at (6.2,3.5) {$T_r$};
\node at (6.1,4.4) {$v_r$};

\draw [line width=1.2pt] (5,5)--(4,4);
\draw [line width=1.2pt] (4,4) circle [radius=.1];
\draw [line width=1.2pt,dashed] (4,4)--(3,3) (4,4)--(5,3);
\draw [line width=1.2pt] (2.8,2.8)--(5.2,2.8)--(5.2,4.2)--(2.8,4.2)--(2.8,2.8);
\node at (4,3.3) {$T_{\ell}$};
\node at (4.1,4.4) {$v_{\ell}$};

\draw [line width=1.2pt] (5,3)--(5.7,2.3);
\draw [fill=red,red] (5.7,2.3) circle [radius=.1];
\draw [red] (5.7,2.3) circle [radius=.2];
\draw [line width=1.2pt,dashed] (5.7,2.3)--(6.4,1.6) (5.7,2.3)--(5,1.6);
\draw [line width=1.2pt] (4.7,2.6)--(6.7,2.6)--(6.7,1.5)--(4.7,1.5)--(4.7,2.6);
\node at (5.9,1.9) {$T_1$};
\node at (6.1,2.3) {$v_1$};

\begin{scope}[shift={(1.4,-1.4)}]
\draw [line width=1.2pt] (5,3)--(5.7,2.3);
\draw [fill=red,red] (5.7,2.3) circle [radius=.1];
\draw [red] (5.7,2.3) circle [radius=.2];
\draw [line width=1.2pt,dashed] (5.7,2.3)--(6.4,1.6) (5.7,2.3)--(5,1.6);
\draw [line width=1.2pt] (4.7,2.6)--(6.7,2.6)--(6.7,1.5)--(4.7,1.5)--(4.7,2.6);
\node at (5.9,1.9) {$T_2$};
\node at (6.1,2.3) {$v_2$};
\end{scope}
\node at (5,-.3) {the decomposition of $T\in\lmark$};
\end{scope}

\begin{scope}[shift={(-8,0)}]
%Branching-labeled trees
\draw [fill=red,red] (15,5) circle [radius=0.1];
\draw [red] (15,5) circle [radius=0.2];
\node at (14.6,5) {$v_1$};

\draw [line width=1.2pt] (15,5)--(16,4);
\draw [fill] (16,4) circle [radius=.1];
\draw [line width=1.2pt,dotted] (16,4)--(16.4,3.6) (16,4)--(15.6,3.6);
\draw [line width=1.2pt] (15.6,4.2)--(16.6,4.2)--(16.6,3.2)--(15.6,3.2)--(15.6,4.2);
\node at (16.1,3.6) {$T_1$};

\draw [line width=1.2pt] (15,5)-- (14,4);
\begin{scope}[shift={(-1,-1)}]
\draw [fill=red,red] (15,5) circle [radius=0.1];
\draw [red] (15,5) circle [radius=0.2];
\draw [line width=1.2pt] (15,5)--(16,4);
\node at (14.6,5) {$v_2$};
\draw [fill] (16,4) circle [radius=.1];
\draw [line width=1.2pt,dotted] (16,4)--(16.4,3.6) (16,4)--(15.6,3.6);
\draw [line width=1.2pt] (15.5,4.2)--(16.5,4.2)--(16.5,3.2)--(15.5,3.2)--(15.5,4.2);
\node at (16,3.6) {$T_2$};
\end{scope}

\draw [line width=1.2pt] (14,4)--(13,3);
\begin{scope}[shift={(-2,-2)}]
\draw [fill=red,red] (15,5) circle [radius=0.1];
\draw [red] (15,5) circle [radius=0.2];
\node at (14.6,5) {$v_3$};
\draw [line width=1.2pt] (15,5)--(16,4);
\draw [fill] (16,4) circle [radius=.1];
\draw [line width=1.2pt,dotted] (16,4)--(16.4,3.6) (16,4)--(15.6,3.6);
\draw [line width=1.2pt] (15.4,4.3)--(16.4,4.3)--(16.4,3.3)--(15.4,3.3)--(15.4,4.3);
\node at (15.9,3.6) {$T_3$};
\end{scope}

\draw [line width=1.2pt,dashed] (13,3)--(12,2);
\begin{scope}[shift={(-3,-3)}]
\draw [line width=1.2pt] (15,5) circle [radius=0.1];
\node at (14.7,5) {$v$};
\draw [line width=1.2pt] (15,5)--(16,4);
\draw [fill] (16,4) circle [radius=.1];
\draw [line width=1.2pt,dotted] (16,4)--(16.4,3.6) (16,4)--(15.6,3.6);
\draw [line width=1.2pt] (15.5,4.2)--(16.5,4.2)--(16.5,3.2)--(15.5,3.2)--(15.5,4.2);
\node at (16,3.5) {$T_r$};
\end{scope}

\draw [line width=1.2pt] (12,2)--(11,1);
\draw [fill] (11,1) circle [radius=.1];
\draw [line width=1.2pt,dotted] (11,1)--(10.6,.6) (11,1)--(11.4,.6);
\draw [line width=1.2pt] (10.5,1.2)--(11.5,1.2)--(11.5,.2)--(10.5,.2)--(10.5,1.2);
\node at (11, .6) {$T_{\ell}$};

\node at (13,-.3) {the decomposition of $T\in\bmark$};
\end{scope}
 \end{tikzpicture}
\caption{Two tree structures that are enumerated by Borel's triangle.}\label{Fig:tree-decomp}
\end{center}
\end{figure}

Eq.\eqref{functional-eq} is useful to determine structures counted by Borel's triangle.  Next  
we describe  tree decompositions for two kinds of marked  binary trees and prove that they are counted by Borel's triangle. 
The first one is the branch-marked binary trees in Theorem \ref{thm:FMS}(b). 
Figure~\ref{Fig:tree-decomp} gives examples of these two decompositions, where the red circle means a mark.

\begin{proof}[Proof of Theorem~\ref{thm:FMS}(b).]
Consider binary trees where only the branching vertices can be marked, denote the set of all such binary tress by $\bmark$. Let $T\in \bmark$ and let $v$ be the first vertex on the left spine of $T$ that is not marked by depth-first search. There exists such a vertex since the last vertex on the left spine is a leave and can not be marked.
Denote by $T_{\ell}$ and $T_r$ the left and right subtrees of $v$ respectively, then they are both branching-labeled binary trees. Now denote the ancestors of $v$ as $v_1, v_2, \ldots$ by depth, every $v_i$ is marked and has a right child. Denote the right subtree of $v_i$ by $T_i$, again we have $T_i\in \bmark$. While $T_{\ell}$ and $T_r$ may be empty, all $T_i$'s must be nonempty, so if we let $\bmm_{n,k}$ be the number of binary trees $T\in\bmark$ with $k$ marked branching vertices and $n$ unmarked vertices, the generating function $B(t,x) = \sum_{n,k\geq 0} \bmm_{n,k}t^k x^n$ also satisfies the relation.
\begin{equation*}
B = 1+x\cdot B^2\cdot \frac{1}{1-t(B-1)}.
\end{equation*}
Together with $b_{1,0}=1$, $b_{2,0}=2$ and $b_{2,1}=1$, we have
 $\bmm_{n,k} = B_{n-1,k}$ and the set  $\bmark$ is enumerated by (the shifted)  Borel's triangle.
\end{proof} 

Our second example is  another  kind of vertex-marked binary trees. 
\begin{thm} \label{Theorem:Borel-marks2} 
The entry $B_{n,k}$ of Borel's triangle counts the set of  binary trees on $n+1$ vertices with $k$ marked vertices such that there is no marked vertices on the right spine. Here the right spine consists of 
 the root and any   vertex $u$  for which all the vertices on the path from the root to $u$ is the right-child of its parent node. 
\end{thm} 
\begin{proof} 
First we consider $\lmark$, the set of all binary trees where the marked vertices are not on its right spine. For any binary tree $T\in \lmark$, let $r$ be the root of $T$ and $v_{\ell}, v_r$ be the left and right children of $r$ respectively, we can find all the marked vertices on the right spine of $v_{\ell}$, name them by depth as $v_1, v_2, \ldots$. In particular, if $v_{\ell}$ is marked, then $v_{\ell} = v_1$. Now remove the edges $(r,v_{\ell}), (r, v_r)$ and all edges $(p(v_i), v_i)$, here $p(v_i)$ denotes the parent of $v_i$, we obtain a forest. Let the tree rooted at $v_r$ in the forest be $T_r$, the tree rooted at $v_{\ell}$ be $T_{\ell}$ and the tree rooted at $v_i$ be $T_i$. While $T_r$ and $T_{\ell}$ might be empty, $T_i$ always contain $v_i$ and thus is not empty. Moreover, $T_r$ and $T_{\ell}$ have no marked vertices on right spine, that is, $T_r, T_{\ell}\in \lmark$. Each $T_i$ has its root $v_i$ marked but no other vertex on the right spine is marked.

Denote by $\lmm_{n,k}$ the number of binary trees $T\in\lmark$ with $k$ marked vertices and $n$ vertices in total.  Hence $\ell_{1,0}=1$, $\ell_{2,0}=2$ and $\ell_{2,1}=1$. 
From the above decomposition, we conclude that the generating function $L(t,x) = \sum_{n,k\geq 0} \lmm_{n,k}t^k x^n$
satisfies the following equation
\begin{equation*}
L = 1+x\cdot L^2\cdot \frac{1}{1-t(L-1)}.
\end{equation*}
which yields   $\lmm_{n,k} = f_{n,k}=B_{n-1,k}$.
\end{proof}

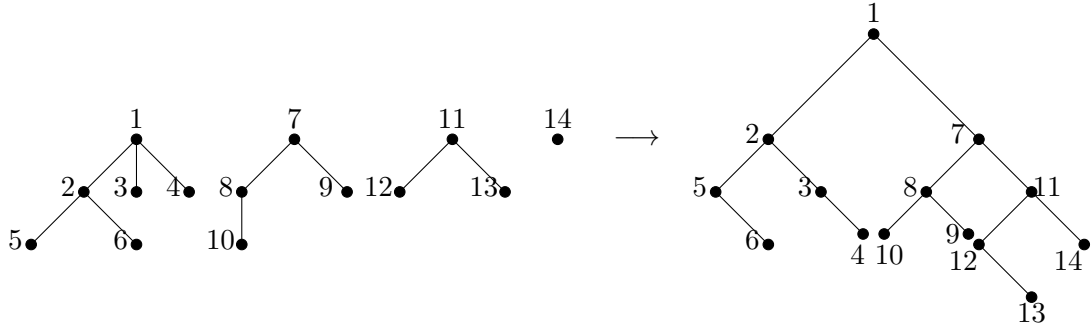
\begin{figure}[htpb]   
\begin{center}
\begin{tikzpicture}[line cap=round,line join=round,>=triangle 45,x=1.0cm,y=1.0cm,scale=.7]
\clip(1.5,1) rectangle (23,8);
%forest
\draw [fill] (4,5) circle [radius=.1] (3,4) circle [radius=.1] (4,4) circle [radius=.1] (5,4) circle [radius=.1];
\draw (4,5)--(3,4) (4,5)--(4,4) (4,5)--(5,4);
\node at (4,5.4) {1};
\node at (2.7,4.1) {2};
\node at (3.7,4.1) {3};
\node at (4.7,4.1) {4};

\draw [fill] (2,3) circle [radius=.1] (4,3) circle [radius=.1];
\draw (3,4)--(2,3) (3,4)--(4,3);
\node at (1.7,3.1) {5};
\node at (3.7,3.1) {6};

\draw [fill] (7,5) circle [radius=.1] (6,4) circle [radius=.1] (8,4) circle [radius=.1] (6,3) circle [radius=.1];
\draw (7,5)--(6,4)--(6,3) (7,5)--(8,4);
\node at (7,5.4) {7};
\node at (5.7,4.1) {8};
\node at (7.6,4.1) {9};
\node at (5.6,3.1) {10};

\draw [fill] (10,5) circle [radius=.1] (9,4) circle [radius=.1] (11,4) circle [radius=.1];
\draw (10,5)--(9,4) (10,5)--(11,4);
\node at (10,5.4) {11};
\node at (8.6,4.1) {12};
\node at (10.6,4.1) {13};

\draw [fill] (12,5) circle [radius=.1];
\node at (12,5.4) {14};

\node at (13.5,5) {$\longrightarrow$};

%binary tree
\draw [fill] (18,7) circle [radius=.1] (16,5) circle [radius=.1] (20,5) circle [radius=.1];
\draw (18,7)--(16,5) (18,7)--(20,5);
\node at (18,7.4) {1};
\node at (15.7,5.1) {2};
\node at (19.6,5.1) {7};

\draw [fill] (15,4) circle [radius=.1] (16,3) circle [radius=.1] (17,4) circle [radius=.1] (17.8,3.2) circle [radius=.1];
\draw (16,5)--(15,4)--(16,3) (16,5)--(17,4)--(17.8,3.2);
\node at (14.7,4.1) {5};
\node at (15.7,3.1) {6};
\node at (16.7,4.1) {3};
\node at (17.7,2.8) {4};

\draw [fill] (19,4) circle [radius=.1] (19.8,3.2) circle [radius=.1] (18.2,3.2) circle [radius=.1];
\draw (20,5)--(19,4)--(19.8,3.2) (19,4)--(18.2,3.2);
\node at (18.7,4.1) {8};
\node at (19.5,3.2) {9};
\node at (18.3,2.8) {10};

\draw [fill] (21,4) circle [radius=.1] (20,3) circle [radius=.1] (21,2) circle [radius=.1] (22,3) circle [radius=.1];
\draw (20,5)--(21,4)--(20,3)--(21,2) (21,4)--(22,3);
\node at (21.3,4.1) {11};
\node at (19.7,2.7) {12};
\node at (21.7,2.7) {14};
\node at (21, 1.7) {13};
\end{tikzpicture}
\end{center}
\caption{The canonical map from plane forests to binary trees}  \label{Fig:tree1}
\end{figure}

\textsc{Remark.} 
The vertex-marked binary trees in Theorem \ref{Theorem:Borel-marks2} can also be obtained from the basic 
marked Catalan structure obtained from item (4) of Theorem \ref{Catalan-array}, which is the set of unlabeled plane forest on $n+1$ vertices with $k$ marked non-root vertices.  To get vertex-marked binary trees,  one uses the canonical bijection from plane forests on $n+1$  vertices to binary trees on $n+1$ vertices, where for any node of a plane forest, its first child becomes the left child, and its first sibling on the right becomes the right child.  Under this bijection  the root of the plane forest becomes the nodes in the right spine the binary tree. 
See Figure~\ref{Fig:tree1} for an illustration.

\end{document}